# Основные конструкции над алгоритмами выпуклой оптимизации и их приложения к получению новых оценок для сильно выпуклых задач

*Гасников А.В. (ПреМоЛаб ФУПМ МФТИ, ИППИ РАН)*
*Камзолов Д.И. (ПреМоЛаб ФУПМ МФТИ)*
*Мендель М.А. (ПреМоЛаб ФУПМ МФТИ)*

**Аннотация**

В статье собраны вместе основные современные конструкции работы с алгоритмами (численными методами) решения задач выпуклой оптимизации. В частности, с помощью искусственного введения неточности в вычисление градиента, следуя Ю.Е. Нестерову, рассматривается "адаптивная игра на гладкости задачи", позволяющая использовать методы настроенные на гладкие задачи для решения негладких задач; рассматривается конструкция рестартов, позволяющая получить из численного метода, ищущего решение задачи выпуклой оптимизации, метод пригодный к использованию для задач сильно выпуклой оптимизации; рассматривается прием регуляризации, позволяющий сводить любую выпуклую задачу к сильно выпуклой. Все эти (и некоторые другие) конструкции (например, композитной оптимизации) описываются исходя из одной общей линии – руководствуясь принципом бритвы Оккама: попытаться изложить современное состояние "оптимальных" численных методов выпуклой оптимизации в пространствах больших размеров (для детерминированных постановок: размерность пространства больше необходимого числа итераций). Статья написана по просьбам коллег и студентов, планирующих использовать собранные в статье конструкции в своей работе.

**Ключевые слова:** композитная оптимизация, быстрый градиентный метод, неточный оракул, универсальный метод Ю.Е. Нестерова, рестарт-техника, регуляризация, mini-batch.

## 1. Введение

В весеннем семестре 2015/2016 учебного года первый автор прочитал курс "Стохастическая и Huge-scale оптимизация" одновременно для студентов Физтеха, Независимого московского университета и студентов магистерской программы ММОС ВШЭ [1]. Также похожие вещи подробно обсуждались с коллегами (А.Ю. Горновым и А.С. Аникиным) в ходе пребывания первого автора в марте 2016 г. в Иркутске. Полученная в ходе этих мероприятий обратная связь привела нас к необходимости записать основные положения, излагаемые в курсе. Несмотря на то, что подавляющее большинство приводимых далее фактов (конструкций) являются хорошо известными (впрочем, при этом большинство приводимых результатов, безусловно, можно назвать современными) в оптимизационном сообществе (прежде всего, благодаря усилиям Ю.Е. Нестерова и А.С. Немировского), мы сочли необходимым собрать часто используемые конструкции в одном месте, и описать их в достаточно популярной форме, удобной для понимания их сути. Мы также постарались разбавить классические результаты некоторыми недавними собственными наработками (многие из этих наработок были получены совместно с П.Е. Двуреченским), которые (на наш взгляд) удачно дополняют имеющиеся результаты, как бы овыпукляя их. В основном новые результаты в данной статье связаны с рассмотрением сильно выпуклых постановок. Отметим также полезные (на наш взгляд) замечания, в которых мы постарались отметить возможные дальнейшие направления развития приводимых в статье результатов.

## 2. Основные результаты

Рассматривается задача выпуклой композитной оптимизации [2]
$$F(x) = f(x) + h(x) \to \min_{x \in Q}. \qquad (1)$$

Мы считаем, что нам доступен $(\delta, L)$-оракул (глава 4 [3], [4]), который, получая на вход произвольный элемент $x \in Q \subseteq \mathbb{R}^n$, выдает такую пару (число + вектор из $\mathbb{R}^n$) $\{f_{\delta,L}(x), g_{\delta,L}(x)\}$, что для любого $y \in Q$

$$0 \le f(y) - f_{\delta,L}(x) - \langle g_{\delta,L}(x), y - x \rangle \le \frac{L}{2}\|y-x\|^2 + \delta.$$

Функция $h(x)$ считается простой структуры (см. [2]), поэтому можно не "обременять" оракула запросами, касающимися этой функции, т.е. мы ее "зашиваем" в шаг метода без обращения к оракулу.

Положим $R^2 = V(x_*, x^0)$, где прокс-расстояние определяется формулой (см., например, главу 2 [3], [4])

$$V(x, z) = d(x) - d(z) - \langle \nabla d(z), x - z \rangle;$$

прокс-функция $d(x) \ge 0$ ($d(x^0) = 0$, $\nabla d(x^0) = 0$) считается сильно выпуклой относительно выбранной нормы $\|\ \|$, с константой сильной выпуклости $\ge 1$; $x_*$ – решение задачи (1) (если решение не единственно, то выбирается то, которое доставляет минимум $V(x_*, x^0)$); $x^0$ – точка старта итерационного процесса; $y^N$ – то, что выдает итерационный процесс после $N$ итераций (обращений к оракулу).

Приводимые далее утверждения фактически получены в работах [3 – 7]. Однако здесь предлагается более наглядная схема рассуждений. Новым в последующем изложении являются результаты, касающиеся сильно выпуклого случая, и следующие после описания этого случая замечания.

**Утверждение 1** (см. главу 4 [3], [4]). *Композитный быстрый градиентный метод (БГМ) Ю.Е. Нестерова с $(\delta, L)$-оракулом (вместо "обычного" оракула, выдающего "настоящие" значения функции и градиенты) сходится ($N$ – число обращений к оракулу) следующим образом (с точностью до констант оценки оптимальны)*

$$F(y^N) - F_* \le \varepsilon, \ N = O\left(\sqrt{\frac{LR^2}{\varepsilon}}\right), \ \delta \le O\left(\frac{\varepsilon}{N}\right). \qquad (2)$$

**Замечание 1** (см. главу 6 [3], [5]). Можно предложить однопараметрическое семейство (параметр $p \in [0,1]$) промежуточных градиентных методов с оценками

$$F(y^N) - F_* \le \varepsilon, \ N = O\left(\left(\frac{LR^2}{\varepsilon}\right)^{\frac{1}{p+1}}\right), \ \delta \le O\left(\frac{\varepsilon}{N^p}\right). \qquad (3)$$

**Утверждение 2** (см. главу 4 [3], [4]). *Пусть*

$$\|\nabla f(y) - \nabla f(x)\|_* \le L_\nu \|y - x\|^\nu \qquad (4)$$

*при некотором* $\nu \in [0,1]$. *Тогда*

$$0 \le f(y) - f(x) - \langle \nabla f(x), y - x \rangle \le \frac{L}{2}\|y - x\|^2 + \delta, \ L = L_\nu \cdot \left[\frac{L_\nu}{2\delta}\frac{1-\nu}{1+\nu}\right]^{\frac{1-\nu}{1+\nu}}. \qquad (5)$$

**Утверждение 3.** *Пусть* $f(x)$ *удовлетворяет условию (4) тогда метод из утверждения 1 с* $(\delta, L)$-*оракулом* $(f(x), \nabla f(x))$ *сходится следующим (оптимальным) образом*

$$F(y^N) - F_* \le \varepsilon, \ N = O\left(\left(\frac{L_\nu R^{1+\nu}}{\varepsilon}\right)^{\frac{2}{1+3\nu}}\right). \qquad (6)$$

*Заметим, что выше в определении* $(\delta, L)$-*оракула* $\delta$ *считается по формуле (2), причем* $N$ *в этой формуле надо подставлять из формулы (6); а* $L$ *считается по формуле (5), в которую надо подставить описанное* $\delta$.

**Схема доказательство.** Воспользуемся утверждениями 1 и 2 (см. формулы (2), (5))

$$N^2 \sim \frac{LR^2}{\varepsilon}, \ L \sim L_\nu \cdot \left(\frac{L_\nu}{\delta}\right)^{\frac{1-\nu}{1+\nu}}, \ \delta \sim \frac{\varepsilon}{N}.$$

$$N^2 \sim \frac{L_\nu^{\frac{2}{1+\nu}} R^2}{\varepsilon \delta^{\frac{1-\nu}{1+\nu}}} \sim \frac{L_\nu^{\frac{2}{1+\nu}} R^2}{\varepsilon^{\frac{2}{1+\nu}} N^{-\frac{1-\nu}{1+\nu}}} \Rightarrow N^{\frac{1+3\nu}{1+\nu}} \sim \frac{L_\nu^{\frac{2}{1+\nu}} R^2}{\varepsilon^{\frac{2}{1+\nu}}} \Rightarrow N \sim \left(\frac{L_\nu R^{1+\nu}}{\varepsilon}\right)^{\frac{2}{1+3\nu}}. \blacksquare$$

**Утверждение 4** (см. описание универсального метода Ю.Е. Нестерова [6]). *Пусть в утверждении 3 используется композитный БГМ с **адаптивным подбором константы Липшица** градиента. Методу известна точность* $\varepsilon > 0$, *с которой хотим решить задачу (больше методу ничего знать о параметрах задачи не обязательно) и доступен оракул, выдающий* $(f(x), \nabla f(x))$. *Тогда оценку (6) можно уточнить следующим образом*

$$F(y^N) - F_* \le \varepsilon, \ N = O\left(\inf_{\nu \in [0,1]} \left(\frac{L_\nu R^{1+\nu}}{\varepsilon}\right)^{\frac{2}{1+3\nu}}\right). \qquad (7)$$

*При этом число обращений за значением функции* $F$ *в среднем будет не более четырех на одну итерацию (итерация характеризуется одним обращение за градиентом).*

**Схема доказательства.** Поясним, что имеется в виду под "композитным БГМ с адаптивным подбором константы Липшица градиента". Подавляющее большинство современных методов содержат в качестве составной части того, что нужно сделать на каждой итерации расчет градиентного отображения (см., например, шаг 2 описанного ниже в п. 3 БГМ – также часто рассматривают проксимальный вариант градиентного отображения, в котором $\|x - x^{k+1}\|^2$ заменяют на $V(x, x^{k+1})$, рассуждения можно перенести и на этот случай). Далее мы приводим процедуру (схожую с оригинальной [6], но все же

отличающуюся от неё выбором $\delta_{k+1}$), которая позволяет "универсализировать" метод за счет изменения шага расчета градиентного отображения.

***База*** $L_0 = 1$.

$(k+1)$***-й шаг.*** Положим $L_{k+1} := L_k/2$. *До тех пор пока*

$$F\left(\operatorname{Grad}_{f,h}^{L_{k+1}}\left(x^{k+1}\right)\right) > \overline{F}_{L_{k+1}}\left(\operatorname{Grad}_{f,h}^{L_{k+1}}\left(x^{k+1}\right); x^{k+1}\right) + \delta_{k+1}, \ \delta_{k+1} \sim \varepsilon^{3/2}/L_{k+1}^{1/2},$$

$$\operatorname{Grad}_{f,h}^{L}\left(x^{k+1}\right) = \arg\min_{x \in Q} \overline{F}_L\left(x; x^{k+1}\right),$$

$$\overline{F}_L\left(x; x^{k+1}\right) = f\left(x^{k+1}\right) + \left\langle \nabla f\left(x^{k+1}\right), x - x^{k+1}\right\rangle + \frac{L}{2}\left\|x - x^{k+1}\right\|^2 + h(x)$$

*выполнять*

$$L_{k+1} := 2L_{k+1}.$$

*Положить*

$$y^{k+1} = \operatorname{Grad}_{f,h}^{L_{k+1}}\left(x^{k+1}\right).$$

Формула для $\delta_{k+1}$ возникла исходя из того, что для используемого базового метода (композитного БГМ) $\delta \sim \varepsilon/N$, а $N \sim \sqrt{L/\varepsilon}$ (см. формулу (2)). На самом деле, надо выбирать $\delta_{k+1}$ следующим (более точным) образом $\delta_{k+1} = \varepsilon\alpha_{k+1}/(2A_{k+1})$, где $A_{k+1} = \sum_{i=1}^{k+1}\alpha_i$, адаптивно подбирая $\alpha_{k+1}$, исходя из соотношения $L_{k+1} = A_{k+1}/\alpha_{k+1}^2$ (см. п. 3). ∎

**Замечание 2.** Если, в свою очередь, вместо обычного оракула, выдающего $(f(x), \nabla f(x))$, использовать $(\delta, L)$-оракул (с константой $L$, определяемой согласно формуле (5) для обобщения (6), и с взятием $\inf_{\nu \in [0,1]}$ в (5) для обобщения (7)), то оценку (7) можно уточнить следующим образом [8] (приведенные оценки не улучшаемы)

$$F\left(y^N\right) - F_* \leq \varepsilon, \ N = \operatorname{O}\left(\inf_{\nu \in [0,1]}\left(\frac{L_\nu R^{1+\nu}}{\varepsilon}\right)^{\frac{2}{1+3\nu}}\right), \ \delta \leq \operatorname{O}\left(\frac{\varepsilon}{N}\right). \tag{8}$$

**Замечание 3.** Исходя из формулы (3) можно провести рассуждения для промежуточного градиентного метода аналогично доказательству утверждения 3. Это приводит в итоге к следующему обобщению формулы (8) [9]

$$F\left(y^N\right) - F_* \leq \varepsilon, \ N = \operatorname{O}\left(\inf_{\nu \in [0,1]}\left(\frac{L_\nu R^{1+\nu}}{\varepsilon}\right)^{\frac{2}{1+2p\nu+\nu}}\right), \ \delta \leq \operatorname{O}\left(\frac{\varepsilon}{N^p}\right), \ p \in [0,1]. \tag{9}$$

Опишем далее основную конструкцию (***рестарт-технику***), позволяющую переносить описанные выше результаты на случай $\mu$-сильно выпуклого функционала $F$ в $\|\ \|$-норме (отметим, что $\mu$-сильно выпуклым может быть только композит $h$).

С учетом утверждения 2 можно ограничиться рассмотрением гладкого случая. Из [2] следует, что БГМ "работает" согласно оценке (левое неравенство имеет место в виду $\mu$-сильно выпуклости $F$)

$$\frac{\mu}{2}\left\|y^{\bar{N}} - x_*\right\|^2 \leq F\left(y^{\bar{N}}\right) - F_* \leq \frac{4LV\left(x_*, x^0\right)}{\bar{N}^2}.$$

Не ограничивая общности [10, 11] будем считать, что прокс-расстояние $V$ можно выбрать так, чтобы оно удовлетворяло условию (в евклидовом случае $\omega_n = 1$)

$$\omega_n = \sup_{x \in Q} \frac{2V(x, x^0)}{\|x - x^0\|^2} = \mathrm{O}(\ln n).$$

Отсюда имеем

$$\|y^{\bar{N}} - x_*\|^2 \le \frac{8LV(x_*, x^0)}{\mu \bar{N}^2} \le \frac{1}{2}\|x^0 - x_*\|^2 \frac{8L\omega_n}{\mu \bar{N}^2}.$$

Выбирая

$$\bar{N} = \sqrt{\frac{8L}{\mu}\omega_n},$$

получим, что

$$\|y^{\bar{N}} - x_*\|^2 \le \frac{1}{2}\|x^0 - x_*\|^2.$$

Выберем в БГМ в качестве точки старта $y^{\bar{N}}$, в качестве прокс-функции (считаем, что так определенная функция корректно определена на $Q$ с сохранением свойства сильной выпуклости)

$$d(x) := d(x - \bar{y}^N + x^0),$$

и снова сделаем $\bar{N}$ итераций, и т.д. Несложно понять, что если мы хотим достичь точности по функции $\varepsilon$, то число $k$ таких рестартов (перезапусков) БГМ достаточно взять (здесь используется стандартное обозначение $\lceil \cdot \rceil$, которое мы поясним примером $\lceil 0.2 \rceil = 1$)

$$k = \left\lceil \log_2\left(\frac{\mu R^2}{\varepsilon}\right) \right\rceil,$$

где в данной формуле $R^2 = \|x^0 - x_*\|^2$ в отличие от формул (2), (3), (6) – (9), в которых $R^2 = V(x_*, x^0)$. Приведенная формула следует из выкладки

$$F\left(\left[y^{\bar{N}}\right]^k\right) - F_* \le \left(\frac{1}{2}\right)^k \|x^0 - x_*\|^2 \frac{4L\omega_n}{\bar{N}^2} = \frac{\mu \|x^0 - x_*\|^2}{2^{k+1}}.$$

Общее число обращений к оракулу будет $N = k\bar{N}$, т.е.

$$N = \sqrt{\frac{8L}{\mu}\omega_n} \left\lceil \log_2\left(\frac{\mu R^2}{\varepsilon}\right) \right\rceil. \tag{10}$$

Эта оценка оптимальна с точностью до множителя $\sim \sqrt{\omega_n}$ [10]. Нам неизвестно, можно ли устранить этот множитель или он неустраним.

**Утверждение 5**. *Используя рестарт-технику (с помощью формулы (10)) можно предложить такое "рестарт обобщение" описанных выше методов (выбирается наиболее общая форма записи (9)), что оценки их работы будут иметь следующий вид*

$$F(y^N) - F_* \le \varepsilon, \quad N = \mathrm{O}\left(\inf_{\nu \in [0,1]} \left(\frac{L_\nu^{\frac{2}{1+\nu}}}{\mu \varepsilon^{\frac{1-\nu}{1+\nu}}} \omega_n\right)^{\frac{1+\nu}{1+2p\nu+\nu}} \cdot \left\lceil \ln\left(\frac{\mu R^2}{\varepsilon}\right) \right\rceil\right), \quad \delta \le \mathrm{O}(\varepsilon/N^p). \tag{11}$$

Данная оценка ранее уже приводилась с $\delta = 0$ при $p = \{0;1\}$ в главе 5 [3]. Здесь мы привели (следуя идее [12]) обобщение соответствующих результатов [3], используя довольно общий прием рестартов (см., например, [13]).

Однако для полноты изложения нам представляется полезным пояснить обратный способ проверки связи формул (9), (11). А именно, далее с помощью конструкции ***регуляризации*** исходной задачи (1) мы опишем общий прием погружения не сильно выпуклых задач в класс сильно выпуклых с последующим извлечением оптимальных методов для не сильно выпуклых задач, исходя из имеющихся оптимальных методов для сильно выпуклых. Кроме того, регуляризация задачи сразу дает эффективный критерий останова метода – в виде контроля малости градиента (градиентного отображения). Не в сильно выпуклом случае такой подход является слишком грубым [14].

Введем семейство $\gamma$-сильно выпуклых в норме $\| \ \|$ задач ($\gamma > 0$)

$$F^{\gamma}(x) = F(x) + \gamma V(x, x^0) \to \min_{x \in Q}. \tag{12}$$

Пусть

$$\gamma \le \frac{\varepsilon}{2V(x_*, x^0)} = \frac{\varepsilon}{2R^2}, \tag{13}$$

и удалось найти $\varepsilon/2$-решение задачи (12), т.е. нашелся такой $y^N \in Q$, что

$$F^{\gamma}(y^N) - F^{\gamma}_* \le \varepsilon/2.$$

Тогда

$$F(y^N) - F_* \le \varepsilon.$$

Действительно,

$$F(y^N) - F_* \le F^{\gamma}(y^N) - F_* \le F^{\gamma}(y^N) - F^{\gamma}_* + \varepsilon/2 \le \varepsilon.$$

Здесь мы использовали определение $F^{\gamma}_*$ и формулу (13)

$$F^{\gamma}_* = \min_{x \in Q}\{F(x) + \gamma V(x, x^0)\} \le F(x_*) + \gamma V(x_*, x^0) \le F_* + \varepsilon/2.$$

Полагая в формуле (11) $\mu = \gamma$, из формулы (13) (понимаемой как равенство), приходим к оценкам формулы (9) с точностью до некоторой степени множителя $\sim \omega_n$:

$$N \sim \left( \frac{L_{\nu}^{\frac{2}{1+\nu}}}{\mu \varepsilon^{\frac{1-\nu}{1+\nu}}} \omega_n \right)^{\frac{1+\nu}{1+2p\nu+\nu}} \sim \left( \frac{L_{\nu}^{\frac{2}{1+\nu}} R^2}{\varepsilon^{\frac{2}{1+\nu}}} \omega_n \right)^{\frac{1+\nu}{1+2p\nu+\nu}} \sim \left( \frac{L_{\nu} R^{1+\nu}}{\varepsilon} \omega_n \right)^{\frac{2}{1+2p\nu+\nu}} \sim \left( \frac{L_{\nu} R^{1+\nu}}{\varepsilon} \right)^{\frac{2}{1+2p\nu+\nu}}.$$

**Замечание 4.** Все описанные методы являются прямо-двойственными (см., например, [15] и цитированную там литературу). В сильно выпуклом случае это не представляет, как правило, особого интереса, потому что наличие сильной выпуклости позволяет эффективно восстанавливать решение соответствующей сопряженной задачи (речь идет о задачах, в которых есть "модель", позволяющая по явным формулам или просто эффективно связывать решения исходной и сопряженной задачи). В не сильно выпуклом случае отмеченное свойство методов представляется весьма полезным.

В виду связи формул (9), (11) может показаться, что не сильно выпуклый случай малоинтересен, поскольку в категориях $O(\ )$ он сводим (с помощью регуляризации) к сильно выпуклому случаю со всеми вытекающими отсюда преимуществами. Отчасти это так, и этим приемом (регуляризацией) часто пользуются. Однако стоит отметить, что процедура рестартов, к сожалению, описана нами таким образом, что на каждом рестарте необходимо делать предписанное число итераций $\bar{N}$ рестартуемого метода. Проблема в том, что в отличие от не сильно выпуклого случая, описанная конструкция рестартов требует явного знания констант Липшица (Гёльдера) градиента гладкой части функционала. В связи с этим могут возникать сложности в практической реализации универсального рестарт-метода, оптимально настраивающегося на параметр $\nu \in [0,1]$. Заметим, что в не сильно выпуклом случае эту задачу решал сам метод. Мы не давали ему на вход никакой информации о гладкости задачи (и тем более, не давали методу константу Липшица (Гельдера) градиента гладкой части функционала), метод сам настраивался на нужную гладкость. Тем не менее, отмеченная проблема разрешима с помощью отслеживания максимальной (на текущем рестарте) константы Липшица градиента гладкой части функционала по накопленной последовательности этих констант. Но не смотря на наличие такого решения, все равно имеются заметные потери в эффективности из-за необходимости выполнения на каждом рестарте предписанного (жестко заданного) числа итераций, которое, как правило, оказывается заметно выше, чем нужно для сокращения расстояния до решения в два раза. На данный момент нам не известны более адаптивные варианты выхода из рестарта, чем выполнения предписанного числа итераций. Естественная тут идея контроля малости градиента (градиентного отображения) приводит в теоретическом плане к заметному ухудшению оценки необходимого числа итераций на каждом рестарте.

В продолжение вышенаписанного заметим, что согласно "теории", если взять не сильно выпуклую задачу, взять не сильно выпуклый метод для этой задачи и применить, то получится оценка времени работы, которая должна совпасть по порядку с оценкой, получаемой при решения регуляризованной задачи с помощью рестартованного метода. На самом деле численные эксперименты говорят о том, что второй подход (при имеющихся сейчас у нас способах программной реализации конструкции рестартов) приводит к увеличению времени работы на один-два порядка, т.е. константа в $O(\ )$ может отличаться на практике в 100 раз и даже больше!

**Замечание 5.** Все описанные методы могут быть распространены (насколько нам известно, это пока еще не сделано в общем случае) на задачи условной оптимизации. Для этого сначала, следуя п. 2.3 [16], стоит рассмотреть минимаксную задачу (ввести правильную лианеризацию исходного функционала и градиентное отображение). Далее использовать идею ***метода нагруженных функционалов*** п. 2.3.4 [16], приводящую к рестартам по неизвестному параметру п. 2.3.5 [16] (оптимальное значение функционала задачи), введение которого, позволяет свести задачу условной оптимизации к минимаксной. Дополнительная плата за такое "введение" (т.е. за рестарты) будет всего лишь логарифмическая, и с точностью до этой "платы" оценки будут оптимальными.

**Замечание 6.** Интересно, на наш взгляд, рассматривать приложения описанных выше методов, где неточность возникает из-за невозможности точного вычисления градиентного отображения (проектирования градиента). В большинстве приложений "стоимость" (время) получения от оракула (роль которого, как правило, играют нами же написанные подпрограммы вычисления градиента) градиента функционала заметно превышает время, затрачиваемое на то, чтобы сделать шаг итерации, исходя из выданного оракулом вектора. Желание сбалансировать это рассогласование (усложнить итерации, сохранив при этом старый порядок сложности, и выиграть за счет этого в сокращении числа итераций), привело к возникновению ***композитной оптимизации*** [2], в которой (аддитивная) часть функционала задачи переносится без лианеризации (запроса градиента) в итерации. Другой способ перенесения части сложности задачи на итерации был описан в замечании 5. Здесь остается еще много степеней свободы, позволяющих играть на том насколько "дорогой" будет оракул и соответствующая (этому оракулу) "процедура проектирования", и том сколько (внешних) итераций потребуется методу для достижения заданной точности. В частности, если обращение к оракулу за градиентом и последующее проектирование требуют, в свою очередь, решения вспомогательных оптимизационных задач, то можно "сыграть" на том, насколько точно надо решать эти вспомогательные задачи, пытаясь найти "золотую середину" между стоимостью итерации и числом итераций. Также можно сыграть и на том, как выделять эти вспомогательные задачи. Другими словами, что понимать под оракулом и под итерацией метода. Общая идея "разделяй и властвуй", применительно к численным методам выпуклой оптимизации может принимать довольно неожиданные и при этом весьма эффективные формы (особенно ярким примером, на наш взгляд, являются методы внутренней точки [11, 16]). Разные варианты описанной игры в связи с транспортно-сетевыми приложениями уже разбирались нами в других работах [8, 14, 17] (см. также п. 3). Интересно было бы систематизировать и структурировать накопленные здесь знания.

**Замечание 7.** Известно, см., например, [10], что для задач стохастической оптимизации все приведенные выше оценки существенным образом модернизируются. И главная особенность этой модернизации заключается в том, что теперь уже гладкость функционала не играет в общем случае существенной роли. Соответствующие оценки (в выпуклом и сильно выпуклом случае) для сходимости в среднем будут иметь вид не зависимо от класса гладкости задачи (см. (4))

$$N = \mathrm{O}\left(\frac{M^2 R^2}{\varepsilon^2}\right), \; N = \mathrm{O}\left(\frac{M^2}{\mu\varepsilon}\omega_n\right). \tag{14}$$

Определение $M = L_0$ см. в формуле (4), что в стохастическом случае можно понимать, как

$$\left\| E_\xi\left[\nabla_x f(y;\xi)\right] - E_\xi\left[\nabla_x f(x;\xi)\right]\right\|_* \le M, \; E_\xi\left[\nabla_x f(x;\xi)\right] = \nabla_x E_\xi\left[f(x;\xi)\right].$$

Оценки достигаются и оптимальны [10] (последняя оценка оптимальна с точностью до множителя $\sim \omega_n$; нам неизвестно, можно ли в принципе избавиться от этого множителя). Естественно задаться вопросом: можно ли улучшить эти оценки, если перейти от выбранных категорий в более точные категории (см. определение $(\delta, L)$-оракула в начале этого пункта): для любого $x \in Q$

$$E_\xi\left[\left\|g_{\delta,L}(x;\xi) - E_\xi\left[g_{\delta,L}(x;\xi)\right]\right\|_*^2\right] \le D,$$

и для любых $x, y \in Q$

$$0 \le f(y) - E_\xi\left[f_{\delta,L}(x;\xi)\right] - \left\langle E_\xi\left[g_{\delta,L}(x;\xi)\right], y - x\right\rangle \le \frac{L}{2}\|y-x\|^2 + \delta?$$

Очевидно, что оценки (14) перестают быть оптимальными в случае, когда $D$ мало. Тут особенно интересны различные приложения, возникающие при специальных рандомизациях, в которых $D$ зависит от $x \in Q$, и по мере приближения к решению $x \to x_*$ выполняется $D(x) \to 0+$ [14, 18–24]. Оказывается, что существует такая вариация описанных ранее методов (описанию планируется посвятить отдельную работу, впрочем, основная идея достаточно простая, и приведена в конце этого замечания), что можно получить следующее обобщение оценки (9) (при $\nu=1$, $p=\{0;1\}$ см. главу 7 [3] и [25], в случае $\nu=1$, $p\in[0,1]$ см. [12])

$$E\left[F\left(y^N\right)\right]-F_* \leq \varepsilon,$$

$$N = \max\left\{\mathrm{O}\left(\underbrace{\inf_{\nu\in[0,1]}\left(\frac{L_\nu R^{1+\nu}}{\varepsilon}\right)^{\frac{2}{1+2p\nu+\nu}}}_{\tilde{N}}\right), \mathrm{O}\left(\frac{DR^2}{\varepsilon^2}\right)\right\}, \quad \delta \leq \mathrm{O}\left(\frac{\varepsilon}{\tilde{N}^p}\right), \quad p\in[0,1].$$

Аналогичное обобщение можно сделать и в сильно выпуклом случае

$$E\left[F\left(y^N\right)\right]-F_* \leq \varepsilon,$$

$$N \simeq \max\left\{\mathrm{O}\left(\underbrace{\inf_{\nu\in[0,1]}\left(\frac{L_\nu^{\frac{2}{1+\nu}}\omega_n}{\mu\varepsilon^{\frac{1-\nu}{1+\nu}}}\right)^{\frac{1+\nu}{1+2p\nu+\nu}}\cdot\left\lceil\ln\left(\frac{\mu R^2}{\varepsilon}\right)\right\rceil}_{\tilde{N}}\right), \mathrm{O}\left(\frac{D}{\mu\varepsilon}\right)\right\}, \quad \delta \leq \mathrm{O}\left(\varepsilon/\tilde{N}^p\right).$$

В приведенных формулах обращает на себя внимание введение $\tilde{N}$. Оказывается с помощью приема, который в западной литературе называют (см., например, [22]) **_mini-batch_** (русского эквивалента пока нет), можно добиться того, чтобы число итераций метода было $\tilde{N}$ [12] (при этом число обращений к оракулу за стохастическим градиентом, по-прежнему, будет $N$). Отсюда и получаются соответствующие (определяемые гладкой частью оценки) оценки на уровень допустимого шума $\delta$ (не случайной природы). Прием заключается в том, что вместо одного раза к оракулу на каждой итерации обращаются $m$ раз за $g_{\delta,L}\left(x;\xi^k\right)$ при одном и том же $x\in Q$, но разных (независимых друг от друга) реализациях $\xi^k$. Исходя из этих обращений на каждой итерации рассчитывается вектор

$$\frac{1}{m}\sum_{k=1}^m g_{\delta,L}\left(x;\xi^k\right),$$

который и используется в качестве "градиента" в методе. Параметр $m\geq 1$ подбирается (минимально возможным), исходя из одного из условий (везде, где мы пишем $\mathrm{O}(\ )$ на самом деле можно писать точные константы $\sim 10^1-10^2$)

$$\tilde{N} \geq \mathrm{O}\left(\frac{DR^2}{m\varepsilon^2}\right), \quad \tilde{N} \geq \mathrm{O}\left(\frac{D}{m\mu\varepsilon}\right).$$

Оказывается, что можно организовать и адаптивный подбор этого параметра (что дает возможность брать в итоговых оценках $\inf$ по $\nu\in[0,1]$) вместе с коэффициентам $\{\alpha_{k+1}\}$ (см. п. 3): $m_{k+1}\simeq 2\alpha_{k+1}D/\varepsilon$.

**Замечание 8.** Отметим, что приводимые выше результаты обобщаются на покомпонентные методы, спуски по направлению и безградиентные методы. То, как преобразуются соответствующие формулы, можно посмотреть, например, в работах [9, 14]. Недавно было обнаружено [20], что все эти обобщения можно получить просто при подстановке в описываемые выше методы (ну а точнее в правильные их модификации), соответствующих рандомизированных вариантов градиента и правильной корректировке способа выбора шагов (с безградиентными методами ситуация чуть посложнее, приходится еще *сглаживать задачу* с помощью свертки с хорошим ядром [10, 20, 21, 24] – не следует путать с методом *двойственного сглаживания* Ю.Е. Нестерова [26]). Здесь имеется некоторый подвох [18], заключающийся в том, что оптимальные оценки для таких рандомизированных вариантов необходимо получать заново (заглядывая в структуру метода), т.е. не пытаться использовать (также оптимальные) оценки, выписанные в замечании 7 (попытка использовать приводит в итоге к заметно завышенным оценкам, не являющимся оптимальными). То есть в данном случае не так уж и просто породить из оптимального метода, скажем, его оптимальный покомпонентный вариант. Имеющиеся сейчас процедуры, которые позволяют это делать (см., например, [18]) требуют аккуратного погружения в структуру метода. Однако при этом, например для (блочно-)покомпонентных методов (в отличие от общих задач стохастической оптимизации), по-видимому, сохраняется возможность перенесения конструкции адаптивного подбора константы Липшица [23] и идеи адаптивной настройки метода на гладкость задачи [18].

### 3. Пример задачи композитной оптимизации (сильно выпуклый случай)

Рассмотрим конкретный пример задачи выпуклой композитной оптимизации [27, 28]

$$F(x) = \frac{1}{2}\|Ax - b\|_2^2 + \mu \sum_{k=1}^n x_k \ln x_k \to \min_{\sum_{k=1}^n x_k = 1,\, x \geq 0}. \qquad (15)$$

Вместо ограничения $\sum_{k=1}^n x_k = 1$ можно рассматривать ограничение $\sum_{k=1}^n x_k \leq 1$.

Разберем два случая а) $0 < \mu \ll \varepsilon/(2\ln n)$ – мало (сильную выпуклость композита в 1-норме можно не учитывать); б) $\mu \gg \varepsilon/(2\ln n)$ – достаточно большое (сильную выпуклость композита в 1-норме необходимо учитывать).

Выберем норму в прямом пространстве $\|\ \| = \|\ \|_1$. Положим

$$f(x) = \frac{1}{2}\|Ax - b\|_2^2,\ h(x) = \mu \sum_{k=1}^n x_k \ln x_k,\ Q = S_n(1) = \left\{x \geq 0 : \sum_{k=1}^n x_k = 1\right\},\ L = \max_{k=1,\ldots,n} \|A^{\langle k \rangle}\|_2^2,$$

где $A^{\langle k \rangle}$ – $k$-й столбец матрицы $A$.

Введем два оператора (см. также утверждение 4)

$$\mathrm{Grad}_{f,h}^L(x^{k+1}) = \arg\min_{x \in Q} \bar{F}(x; x^{k+1}),$$

$$\bar{F}(x; x^{k+1}) = f(x^{k+1}) + \langle \nabla f(x^{k+1}), x - x^{k+1}\rangle + \frac{L}{2}\|x - x^{k+1}\|_1^2 + h(x);$$

$$\mathrm{Mirr}_{f,h,z^k}^\alpha(\nabla f(x^{k+1})) = \arg\min_{x \in Q}\left\{\langle \nabla f(x^{k+1}), x - z^k\rangle + \frac{1}{\alpha}V(x, z^k) + h(x)\right\},$$

где прокс-расстояние (расстояние Брэгмана) определяется формулой [11, 26, 29]

$$V(x, z) = d(x) - d(z) - \langle \nabla d(z), x - z\rangle,$$

прокс-функция $d(x) \geq 0$ считается сильно выпуклой относительно выбранной нормы $\|\ \| = \|\ \|_1$, с константой сильной выпуклости $\geq 1$. Для случая а) можно выбирать

$$d(x) = \ln n + \sum_{k=1}^{n} x_k \ln x_k.$$

Тогда

$$V(x,z) = \sum_{k=1}^{n} x_k \ln(x_k/z_k), \ R^2 \leq \ln n.$$

Существуют варианты БГМ (см., например, алгоритм 8 главы 2 [3]), в которых шаг $\text{Grad}_{f,h}^L(x^{k+1})$ заменяется его проксимальным аналогом, т.е., грубо говоря, $\|x - x^{k+1}\|_1^2$ заменяется в выражении для $\bar{F}(x; x^{k+1})$ на $V(x, x^{k+1})$. В таком варианте метода (оценки скорости сходимости аналогичны оценкам, приводимым в утверждении 6) мы имеем ситуацию, когда композит совпадает по форме с прокс-расстоянием (энтропийного типа), и шаг итерации осуществим по явным формулам (см., например, [11, 30]). Таким образом, стоимость итерации будет $\text{O}(nnz(A))$, где $nnz(A)$ – число ненулевых элементов в матрице $A$ (считаем, что это число $\geq n$).

Для случая б) планируется использовать рестарт-технику (см. п. 2). Но для выбранной функции $V(x,z)$ (расстояние Кульбака–Лейблера [30]) процедура рестартов некорректна. Однако существует другой способ выбора прокс-функции (детали, см., например, в работах [10, 11, 13, 30])

$$d(x) = \frac{1}{2(a-1)} \|x\|_a^2, \ a = \frac{2 \ln n}{2 \ln n - 1}. \tag{16}$$

В этом случае $R^2 = \text{O}(\ln n)$, $\omega_n = \text{O}(\ln n)$.

Опишем вариант быстрого градиентного метода Ю.Е. Нестерова в форме [29]. Здесь мы распространяем подход работы [29] на задачи композитной оптимизации. Фактически предложенный далее алгоритм есть сочетание БГМ работы [29] и конструкции композитной оптимизации работы [2].

Определим две числовые последовательности $\{\alpha_k, \tau_k\}$:

$$\alpha_1 = \frac{1}{L}, \ \alpha_{k+1} = \frac{1}{2L} + \sqrt{\frac{1}{4L^2} + \alpha_k^2}, \ \tau_k = \frac{1}{\alpha_{k+1} L}.$$

В случае адаптивного подбора константы Липшица (см. утверждение 4)

$$\alpha_1 = \frac{1}{L_1}, \ \alpha_{k+1} = \frac{1}{2L_{k+1}} + \sqrt{\frac{1}{4L_{k+1}^2} + \alpha_k^2 \frac{L_k}{L_{k+1}}}, \ \tau_k = \frac{1}{\alpha_{k+1} L_{k+1}}.$$

Заметим, что при $k \gg 1$

$$\alpha_k \sim \frac{k}{2L}, \ \tau_k \sim \frac{2}{k}.$$

БГМ

1. $x^{k+1} = \tau_k z^k + (1 - \tau_k) y^k$.
2. $y^{k+1} = \text{Grad}_{f,h}^L(x^{k+1})$.

3. $z^{k+1} = \text{Mirr}_{f,h,z^k}^{\alpha_{k+1}}\left(\nabla f\left(x^{k+1}\right)\right)$.

4. Если не выполняется критерий останова (можно по-разному определять [14]), положить
$$k := k+1$$
и перейти к п.1. Иначе остановиться и выдать $y^{k+1}$.

**Утверждение 6.** *Для задачи (15) БГМ генерирует такую последовательность точек $\left\{x^k, y^k, z^k\right\}_{k=0}^{N}$, что имеют место следующие неравенства (второе неравенство означает, что описанный вариант БГМ является прямо-двойственным методом)*

$$F\left(y^N\right) - F_* \leq \frac{4LR^2}{(N+1)^2}.$$

$$\alpha_N^2 L F\left(y^N\right) \leq \min_{x \in Q} \left\{ \sum_{k=0}^{N-1} \alpha_{k+1} \left\{ f\left(x^{k+1}\right) + \left\langle \nabla f\left(x^{k+1}\right), x - x^{k+1} \right\rangle + h(x) \right\} + V\left(x, x^0\right) \right\}.$$

**Замечание 9.** Заметим, что если условие $x \in Q$ можно записать, например, как $g(x) \leq 0$, и ввести двойственную функцию

$$G(\lambda) = \min_x \left\{ f(x) + h(x) + \left\langle \lambda, g(x) \right\rangle \right\}, \; \lambda \geq 0,$$

то, поскольку

$$\alpha_N^2 L = \sum_{k=0}^{N-1} \alpha_{k+1} \stackrel{def}{=} S_N \text{ и } f\left(x^{k+1}\right) + \left\langle \nabla f\left(x^{k+1}\right), x - x^{k+1} \right\rangle \leq f(x),$$

получим

$$F\left(y^N\right) - G\left(\bar{\lambda}^N\right) \leq \frac{4L\tilde{R}^2}{(N+1)^2},$$

где $\bar{\lambda}^N = \lambda^N / S_N$, $\lambda^N$ – множитель Лагранжа к ограничению $g(x) \leq 0$ в задаче

$$\sum_{k=0}^{N-1} \alpha_{k+1} \left\{ f\left(x^{k+1}\right) + \left\langle \nabla f\left(x^{k+1}\right), x - x^{k+1} \right\rangle + h(x) \right\} + V\left(x, x^0\right) \to \min_{g(x) \leq 0},$$

а $\tilde{R}^2 = V\left(\tilde{x}, x^0\right)$, где $\tilde{x}$ – решение задачи

$$f(x) + h(x) + \left\langle \bar{\lambda}^N, g(x) \right\rangle \to \min_x.$$

Последнее условие, к сожалению, в общем случае не дает возможности как-то разумно оценивать сверху $\tilde{R}^2$, как следствие, возникает проблема с теоретической оценкой зазора двойственности. Проблема решается, если удается обосновать возможность компактификации. Пример того, как эту компактификацию можно делать (на основе "Слейтеровских соображений") будет описан далее (см. формулы (20), (21)).

Используя описанную в предыдущем пункте технику рестартов можно получить из утверждения 6 его аналог в случае сильно выпуклой постановки задачи (15) – случай б). Мы опускаем соответствующие рассуждения, и остановимся подробнее на том, как осуществлять шаг итерации описанного БГМ в случае б), т.е. когда прокс-функция выбирается согласно (16). Сложность выполнения одной итерации (дополнительная к вычислению градиента гладкой части функционала $\text{O}(nnz(A))$) определяется тем, насколько эффективно можно решить задачу следующего вида

$$\tilde{F}(x) = \langle c, x \rangle + \|x\|_a^2 + \bar{\mu} \sum_{k=1}^{n} x_k \ln x_k \to \min_{x \in S_n(1)}. \quad (17)$$

Задачу (17) удобно переписать следующим почти "сепарабельным" образом

$$\langle c, x \rangle + t + \bar{\mu} \sum_{k=1}^{n} x_k \ln x_k \to \min_{\substack{x \in S_n(1), \|x\|_a^a \leq t^{a/2}, \\ 0 \leq t \leq n^{2/a},\, 0 \leq x_k \leq 1,\, k=1,\ldots,n}}$$

Слово "почти" можно убрать, если с помощью метода множителей Лагранжа переписать задачу следующим образом

$$\tilde{G}(\lambda) = \min_{\substack{0 \leq t \leq n^{2/a}, \\ 0 \leq x_k \leq 1,\, k=1,\ldots,n}} \left\{ \sum_{k=1}^{n} c_k x_k + t + \lambda_1 \cdot \left( \sum_{k=1}^{n} x_k - 1 \right) + \lambda_2 \cdot \left( \sum_{k=1}^{n} x_k^a - t^{a/2} \right) + \bar{\mu} \sum_{k=1}^{n} x_k \ln x_k \right\} \to \max_{\lambda_1 \in \mathbb{R},\, \lambda_2 \geq 0}. \quad (18)$$

Поиск минимума $(x(\lambda), t(\lambda))$, где

$$t(\lambda) = \min\left\{ \left( \frac{\lambda_2 a}{2} \right)^{\frac{2}{2-a}}, n^{\frac{2}{a}} \right\},$$

сводится к решению $n$ одномерных задач сильно выпуклой оптимизации на отрезке $[0,1]$. Таким образом, если задаться некоторой точностью $\sigma > 0$, то за время $\mathrm{O}(n \ln(n/\sigma))$ методом деления отрезка пополам (или, скажем, методом золотого сечения [31]) можно найти такой $\tilde{x}^\sigma(\lambda)$, что

$$\|\tilde{x}^\sigma(\lambda) - x(\lambda)\|_1 = \mathrm{O}(\sigma). \quad (19)$$

Далее попробуем (следуя [32]) оценить "запас" в **условии Слейтера**, чтобы, исходя из этого, оценить сверху размер решения $\lambda = \lambda_*$ двойственной задачи (18) (в приводимой далее выкладке, приводящей к формуле (20), для упрощения записи мы опускаем нижний индекс "*" у $\lambda$). Из сильной двойственности имеем

$$-\|c\|_\infty - \bar{\mu} \ln n \leq \tilde{F}_* = \tilde{G}_* \leq \sum_{k=1}^{n} c_k \bar{x}_k + \bar{t} + \lambda_1 \cdot \left( \sum_{k=1}^{n} \bar{x}_k - 1 \right) + \lambda_2 \cdot \left( \sum_{k=1}^{n} \bar{x}_k^a - \bar{t}^{a/2} \right) + \bar{\mu} \sum_{k=1}^{n} \bar{x}_k \ln \bar{x}_k.$$

Если $\lambda_1 \geq 0$, то положим $\bar{t} = 1$, $\bar{x}_k = 1/(2n)$, $k = 1,\ldots,n$. Тогда

$$\frac{1}{2}\lambda_1 + \frac{1}{2}\lambda_2 \leq 2\|c\|_\infty + 2\bar{\mu}\ln(n) + 1,$$

Если $\lambda_1 < 0$, то положим $\bar{t} = 8$, $\bar{x}_k = 2/n$, $k = 1,\ldots,n$. Тогда

$$|\lambda_1| + \lambda_2 \leq 3\|c\|_\infty + 2\bar{\mu}\ln(2n) + 8.$$

В любом случае, с хорошим запасом можно гарантировать, что

$$\|\lambda_*\|_1 \leq 4\|c\|_\infty + 4\bar{\mu}\ln(2n) + 8 \overset{def}{=} C. \quad (20)$$

Таким образом, чтобы решить задачу (17), мы должны решить двойственную задачу (18), которую (в виду формулы (20)) можно переписать следующим образом

$$\breve{G}(\lambda) = -\tilde{G}(\lambda) \to \min_{\substack{\lambda_1 \in \mathbb{R},\, \lambda_2 \geq 0 \\ \|\lambda\|_1 \leq C}}. \quad (21)$$

Поскольку эта задача оптимизации на двумерной плоскости (т.е. в пространстве малой размерности), то ее можно решать, скажем, методом эллипсоидов [10]. При этом для расчета градиента $\breve{G}(\lambda)$ мы должны решить задачу (17) и воспользоваться формулой Демьянова–Данскина [32]

$$\frac{\partial \breve{G}}{\partial \lambda_1} = 1 - \sum_{k=1}^{n} x_k(\lambda), \quad \frac{\partial \breve{G}}{\partial \lambda_2} = t(\lambda)^{a/2} - \sum_{k=1}^{n} x_k(\lambda)^a.$$

К сожалению, точно решить задачу (18) мы не можем, зато можем найти приближенное значение градиента. Точнее говоря, в виду (19), (20), мы можем найти для задачи (21) $\delta = \mathrm{O}(C\sigma)$-градиент $\nabla_\delta G(\lambda)$ (см., например, [33]). Если использовать в методе эллипсоидов в пространстве размерности $r$ (в нашем случае $r = 2$) вместо градиента $\delta$-градиент (чаще говорят $\delta$-субградиент, но в нашем случае можно говорить о градиенте), то имеют место следующие оценки [10]

$$\breve{G}(\lambda^N) - \breve{G}_* \leq \varepsilon, \quad N = \mathrm{O}(r^2 \ln(C/\varepsilon)), \quad \delta \leq \mathrm{O}(\varepsilon), \tag{22}$$

При этом стоимость одной итерации будет $\mathrm{O}(r^2)$. Число итераций можно сократить в $\sim r$ раз, сохранив сложность итерации (см., например, [22]). В нашем случае стоимость одной итерации будет $\mathrm{O}(n \ln(nC/\varepsilon))$.

Однако решение задачи (18) (или (21)), в смысле (22) еще не гарантирует возможность точного восстановления решения задачи (17). Для того чтобы показать, что метод эллипсоидов с той же по порядку точностью $\varepsilon$ позволяет восстанавливать (без каких бы то ни было существенных дополнительных затрат) решение задачи (17) нужно воспользоваться **прямо-двойственностью** этого метода [34]. В виду компактности множества (единичный симплекс), на котором ведется оптимизация в прямом пространстве и сильной выпуклости функционала прямой задачи (17) мы не просто восстанавливаем из прямо-двойственной процедуры метода эллипсоидов решение задачи (17) с точностью по функционалу (прямой задачи) порядка $\varepsilon$, но и делаем это в нужном нам более сильном смысле – см. п. 5.5.1 (следует сравнить с п. 4.6 [3] и п. 2 выше). Формула 5.5.15 [11] гарантирует при этом справедливость следующего результата.

**Утверждение 7.** *Для задачи (15) в случае б) БГМ с рестартами и с прокс-функцией (16) приводит к необходимости на каждой итерации наряду с расчетом градиента гладкой части функционала ($\mathrm{O}(nnz(A))$ операций) два раза решать задачу типа (17) с помощью перехода к двойственной задаче и ее решения с помощью прямо-двойственной версии метода эллипсоидов ($\mathrm{O}(n \ln(C/\varepsilon) \ln(nC/\varepsilon))$ операций). При этом*

$$F(y^N) - F_* \leq \varepsilon,$$

*если общее число итераций (обращений к оракулу за градиентом)*

$$N = \mathrm{O}\left(\sqrt{\frac{L}{\mu} \ln n} \left\lceil \ln\left(\frac{\mu}{\varepsilon}\right) \right\rceil\right) = \mathrm{O}\left(\sqrt{\frac{\max_{k=1,\ldots,n} \|A^{\langle k \rangle}\|_2^2 \ln n}{\mu}} \left\lceil \ln\left(\frac{\mu}{\varepsilon}\right) \right\rceil\right). \tag{23}$$

Заметим, что в "пороговой" ситуации, отвечающей регуляризации (см. п. 2), $\mu \simeq \varepsilon/(2\ln n)$. В этом случае формула (23) примет вид

$$N = \mathrm{O}\left(\sqrt{\frac{\max\limits_{k=1,\ldots,n}\left\|A^{\langle k\rangle}\right\|_2^2 \ln^2 n}{\varepsilon}}\right),$$

что с точностью до $\sim \sqrt{\ln n}$ соответствует оценке в случае а). Отличие случая а) и б) также и в том, что в случае а) существует способ добиться стоимости итерации $\mathrm{O}(nnz(A))$, а в случае б) нам не известно более эффективного способа, чем способ (описанный выше) со стоимостью итерации

$$\mathrm{O}\bigl(nnz(A) + n\ln(C/\varepsilon)\ln(nC/\varepsilon)\bigr).$$

Поскольку в типичных приложениях первое слагаемое заметно доминирует второе, то можно было бы не сильно задумываться (что часто и делают на практике) о плате за невозможность выполнения "проектирования" по явным формулам и не сильно задумываться с какой точностью решать вспомогательную задачу, делая это с точностью длины мантиссы (описанный подход позволяет решать ее с очень хорошей точностью, и сложность решения вспомогательной задачи практически не чувствительна к этой точности). По-видимому, этот тезис имеет достаточно широкий спектр практических приложений. Настоящий пункт имел одной из своих целей на конкретном примере более подробно, чем это принято на практике, продемонстрировать, тезис, что для большого класса задач наличие явных формул для шага итерации не есть сколько-нибудь сдерживающие обстоятельство для использования метода. Используемая при этом техника и способ рассуждений характерным образом (на наш взгляд) демонстрируют современный арсенал средств (описанных в п. 2) решения задач выпуклой оптимизации в пространствах больших размеров.

### 4. Заключение

Когда статья готовилась к печати, авторам стала известна работа [35], в которой предложено изящное обобщение конструкций рестартов/регуляризации и сглаживания.

В целом данная работа писалась не столько как научная статья (хотя она и содержит новые результаты, прежде всего, в утверждениях 5, 7 и ряде замечаний), а, скорее, как удобный для использования (хотелось бы на это надеяться) материал обзорного характера, содержащий основные конструкции (или хотя бы ссылки на них) современных численных методов выпуклой оптимизации в пространствах больших и сверх больших размеров. Мы специально довольно много места отвели описательной части (содержащей мало строгих фактов, но много идей), в которой постарались обрисовать возможные направления развития соответствующих методов и подходов. Мы надеемся, что это побудит читателей к работе в указанных направлениях.

При отборе материала мы исходили из желания продемонстрировать целостную картину с проработкой различных связей. Тем не менее, большое количество деталей, к сожалению, мы вынуждены были опустить, стараясь приводить в соответствующих местах наиболее удобные для восстановления результата ссылки.

**Литература**